\documentclass[doublespace]{amsart}

\usepackage[english]{babel}
\usepackage[T1]{fontenc}
\usepackage[utf8]{inputenc}
\usepackage{amsmath,amssymb,amsthm}
\usepackage{textcomp}
\usepackage{stmaryrd}
\usepackage[sc]{mathpazo}
\usepackage[cal=boondox]{mathalfa}
\usepackage{euler}

\usepackage{mathtools}
\usepackage{comment}
\usepackage{enumerate}

\usepackage[all]{xy}
\usepackage{hyperref}
\usepackage{tikz}
\usetikzlibrary{cd}
\usetikzlibrary{matrix}
\usetikzlibrary{knots,hobby,decorations.pathreplacing,shapes.geometric,calc,decorations.text,shapes.misc,calc,math,arrows,decorations.markings}

\definecolor{couleur3}{rgb}{0.0, 0.5, 1.0} % bleu
\definecolor{couleur4}{rgb}{0.55, 0.71, 0.0} % vert
\definecolor{couleur1}{rgb}{0.87, 0.45, 1.0} % rose
\definecolor{couleur2}{rgb}{0.81, 0.06, 0.13} % rouge

\usepackage{pgfplots}
\pgfplotsset{compat=1.15}
\usepackage{mathrsfs}

\usepackage[nameinlink]{cleveref}

\usepackage{todonotes}
\setuptodonotes{textcolor=blue, linecolor=blue, bordercolor=blue, backgroundcolor=white, size=footnotesize}
\usetikzlibrary{arrows}

\renewcommand{\k}{k}

\newcommand{\G}{\mathbb{G}}

\newcommand{\N}{\mathbb{N}}

\renewcommand{\P}{\mathbb{P}}
\newcommand{\Q}{\mathbb{Q}}
\newcommand{\R}{\mathbb{R}}

\newcommand{\Z}{\mathbb{Z}}

\renewcommand{\hbar}{\overline{h}}

\DeclareMathOperator{\Hom}{\mathrm{Hom}}

\DeclareMathOperator{\Spec}{\mathrm{Spec}}
\DeclareMathOperator{\Span}{\mathrm{Span}}

\DeclareMathOperator{\Cone}{\mathrm{Cone}}

\DeclareMathOperator{\Int}{Int}

\newcommand{\inj}{\hookrightarrow}

\crefname{subsection}{subsection}{the subsections}
\crefname{section}{section}{the sections}
\crefname{subsubsection}{subsubsection}{the subsubsections}

\crefname{Thm}{theorem}{theorems}
\Crefname{Thm}{Theorem}{Theorems}
\crefname{Prop}{proposition}{propositions}
\Crefname{Prop}{Proposition}{Propositions}
\crefname{Lemme}{lemma}{lemmas}
\Crefname{Lemme}{Lemma}{Lemmas}
\crefname{Cor}{corollary}{corollaries}
\Crefname{Cor}{Corollary}{Corollaries}
\crefname{Const}{construction}{constructions}
\Crefname{Const}{Construction}{Constructions}
\crefname{Ex}{example}{examples}
\Crefname{Ex}{Example}{Examples}
\crefname{Def}{definition}{definitions}
\Crefname{Def}{Definition}{Definitions}
\crefname{Not}{notation}{notations}
\Crefname{Not}{Notation}{Notations}
\crefname{Rem}{remark}{remarks}
\Crefname{Rem}{Remarque}{Remarques}

\begin{document}
\title{Seminormal toric varieties}
\author{François BERNARD and Antoine BOIVIN}
\email{\href{mailto:antoine.boivin@univ-angers.fr}{antoine.boivin@univ-angers.fr}}
\email{\href{mailto:fbernard@imj-prg.fr}{fbernard@imj-prg.fr}}

\date\today


\begin{abstract}
    In this paper, we provide a combinatorial description of seminormal toric varieties. The corresponding combinatorial object is a fan equipped with a collection of groups assigned to each cone. This framework introduces a more general class of toric varieties than classical normal toric varieties, while having simpler combinatorial data compared to general non-normal toric varieties.
\end{abstract}

\maketitle

\newtheorem{Thm}{Theorem}[subsection]
\newtheorem{Prop}[Thm]{Proposition}
\newtheorem{Lemme}[Thm]{Lemma}
\newtheorem{Cor}[Thm]{Corollary}
\newtheorem{Const}[Thm]{Construction}
\newtheorem*{rappel}{Theorem}

\theoremstyle{definition}

\newtheorem{Ex}[Thm]{Example}
\newtheorem{Def}[Thm]{Definition}
\newtheorem{Not}[Thm]{Notation}

\theoremstyle{remark}
\newtheorem{Rem}[Thm]{Remark}

\newcommand{\fait}[1]{

\todo[bordercolor = black, textcolor =white, color = orange]{Fait : #1}

}

\setcounter{secnumdepth}{4}

\tableofcontents

\section*{Introduction}

Toric varieties are a very useful and classical tool in algebraic geometry because they can be described through combinatorics with fans. In order to have this correspondence, non-affine toric varieties are almost exclusively assumed to be normal. However, one might wonder whether some combinatorial description still holds if we relax the normality condition. For instance, non-normal toric varieties can sometimes arise as flat degenerations of algebraic varieties (see, for example, \cite{Anderson2013ToricDegen}, \cite{Harada2015ToricDegen}). In this paper, we describe the combinatorics of seminormal toric varieties.

The concept of seminormalization was first introduced in the context of moduli problems by Andreotti and Norguet \cite{AndreottiNorguet} for analytic varieties, later extended by Traverso \cite{TraversoPicard} to general rings in the study of Picard groups, and by Andreotti and Bombieri \cite{AndreottiBombieri} for schemes. It was further explored by Leahy and Vitulli \cite{LeahyVitulli, VitulliCorrections87} and by Greco and Traverso \cite{GrecoTraverso} in the context of algebraic varieties. For a variety $X$ defined over an algebraically closed field of characteristic $0$, its seminormalization $X^+$ is the largest intermediate variety between $X$ and its normalization that remains bijective to $X$. In this setup, it has been studied in \cite{Bernard2021} by using rational functions that extend continuously, for the strong topology, to the closed points of the variety, and it is used in \cite{BFMQ} to characterize homeomorphisms between algebraic varieties. The seminormalization also appears in the study of singularities in the minimal model program (see \cite{KollarSingularMMP} and \cite{KollarVariantsOfNormality}). It has the property of having "multicross" singularities in codimension 1 (see \cite{Multicross}), meaning that they are locally analytically isomorphic to the union of linear subspaces of affine space meeting transversally along a common linear subspace. Seminormalization also finds recent applications in real algebraic geometry \cite{SeminormRealAlgGeo, MonnierCentralGeo, BernardRealSN} where it is related to the theory of "regulous" functions introduced in \cite{KollarNowak2015ContinuousRational,ArticleRegulu}.

In commutative algebra, several definitions of seminormalization coexist. Although geometrically intuitive, the definition of Traverso has been conveniently replaced in \cite{Hamann75, LeahyVitulli, SwanOnSeminormality} by some arithmetical definitions which are algebraically more convenient to use and which are all equivalent to the one of Traverso. In particular, a ring $A$ is seminormal in $B$ if every element $b\in B$ verifying $b^2, b^3\in B$, belongs to $A$. This is equivalent to asking that, for every element $b\in B$, if there exists $m\in \N$ such that $b^n\in B$ for all $n>m$, then $b$ belongs to $A$. We refer to Vitulli \cite{VitulliSurvey2011} for an extensive survey from 2011 on seminormality and to the Macaulay2 package \cite{PackageMacauley2Seminormal} to do explicit computations of the seminormalization of reduced rings.

Considering these definitions, the notion of seminormal monoids follows naturally: $S$ is seminormal in $\text{gp}(S)$ if, for all $m\in \mathrm{gp}(S)$, then $2m,3m \in S$ implies $m\in S$. So that a monoid algebra $k[S]$ is seminormal if and only if $S$ is seminormal. Seminormal monoids were first considered by Gubeladze \cite{Gubeladze1989AndersonsConj} in order to prove a conjecture of Anderson \cite{Anderson1978ProjectiveModules} asking for which class of monoids, all the projective modules over the corresponding monoid algebras are free. Seminormal monoids were further studied by Swan in \cite{Swan1992Gubeladze}, and by Reid and Roberts in \cite{REID2001703} who gave a geometric description of the seminormalization of a monoid. Namely, it consists in saturating the interior of each face of the monoid. Conversely, in the case of affine monoids (i.e., affine toric varieties), any seminormal monoid can be obtained by associating a group to each of its faces. In this paper, we extend this description to the non-affine case by associating a group to each cone of a given fan. To achieve this, we use a combinatorial construction introduced by Gonzales Perez and Teissier \cite{TeissierGonzales2012}, known as \textit{fans with attached monoids} which correspond to toric varieties admitting a good action. This construction consists, roughly speaking, of a fan along with a collection of monoids that satisfy gluing conditions associated to the cones of the fan. In the case of seminormal toric varieties, the attached monoids are seminormal, which simplifies the combinatorial data and the gluing conditions. The goal of this paper is to describe precisely these combinatorial data, thereby providing a better understanding of seminormality and introducing a more general class of toric varieties than the classical normal toric varieties, which have simpler combinatorial data compared to general non-normal toric varieties with good action.\\

The article is organized as follows. In \Cref{SectionSemisaturatedMonoids} we present seminormal monoids and the theorem of Reid and Roberts. We choose to call them "semisaturated" monoids since normal toric varieties correspond to "saturated" monoids. Then we prove several necessary lemmas, which, in particular, extend remark 4.4 of \cite{REID2001703} and show that a seminormal affine monoid is completely characterized by a collection of groups associated to each of its faces. In \Cref{SectionSeminormalAffTorVar} we recall the role of normality in the theory of toric geometry and we present seminormal affine toric varieties. In \Cref{SectionGeneralTorVar} we prove the main theorem of the paper. We recall briefly, in \Cref{SubsecTeissierGonzPerez}, the construction of fans with associated monoids of Gonzales Perez and Teissier. Finally, in \Cref{SubsecFansWithAttachedGroups}, we introduce the notion of fans with attached groups which consists of a triple $(N,\Sigma,G)$ where $N$ is a lattice, $\Sigma$ is a fan of $N_{\R}$, and $G$ is a collection of groups $\{G_{\sigma}\}_{\sigma\in\Sigma}$ such that $G_0$ is the dual lattice of $N$ and, for every $\tau\prec\sigma$, the group $G_{\sigma}$ is a finite index subgroup of $G_{\tau}\cap\sigma^\perp$. Then we define the notion of morphisms of fans with attached groups and we show that there is a well-defined functor between this category and the category of fans with attached semisaturated monoids. For a fan $(N,\Sigma,G)$ with attached groups, the functor associates to each face $\sigma\in \Sigma$ a semisaturated monoid $\Gamma_\sigma$ as follows
\[
    \Gamma_\sigma\coloneqq\coprod_{\tau \preceq \sigma} G_\tau \cap \Int(\sigma^\vee \cap \tau^\perp).
\]
By showing that this functor is an equivalence of categories, we get the main theorem of this paper:
\begin{rappel}
    There is an equivalence of categories between the category $\mathfrak{TorVar}^{sn}$ of seminormal toric varieties with a good action and the category $\mathfrak{Fan}^{gr}$ of fans with attached groups.
\end{rappel}

In all the paper, the varieties are taken over any field $\k$, the cones are supposed rational, the monoids will be commutative, cancellative and finitely generated (i.e. monoids $S$ which can be embedded in a finitely generated abelian group).
 
If $S \subset M$ is a set included in a group, we will note $\N S$ (resp. $\Z S$) the abelian monoid (resp. abelian group) generated by the elements of $S$ in $M$ and if $S \subset V$ is a set included in a $\R$-vector space $V$, we will note $\R_{\geq 0} S$ or $\Cone(S)$ the cone generated by the elements of $S$ in $V$ i.e. the cone of finite non-negative linear combinations of elements of $S$. Finally, if $S$ is a subset of a group, we note $k[S]$ the monoid $\k$-algebra of $\N S$ over $k$.\\

\textbf{Acknowledgement :} The first author was supported by Plan d’investissements France 2030, IDEX UP ANR-18-IDEX-0001 and the second author was supported by the Région Pays de la Loire from the France 2030 program, Centre Henri Lebesgue ANR-11-LABX-0020-01.

\section{Semisaturated monoids}\label{SectionSemisaturatedMonoids}

We begin by presenting semisaturated monoids and the Reid-Roberts theorem, which states that a monoid is semisaturated if the interior of each face corresponds to the interior of its saturation. We then build upon a remark from \cite{REID2001703}, that semisaturated monoids can be constructed from a collection of groups. Our goal in developing this idea is to extend later this construction to the non-affine case.

\begin{Lemme}[\cite{REID2001703} Theorem 3.2 b)]\label{LemmaSemisaturatedEquivalences}

    Let $S$ be a (cancellative) monoid and $M$ be a group such that $S \subset M$. Then the following properties are equivalent:%\todo{cancellative monoid ?}
    \begin{enumerate}
    %\item $ \forall m\in M, \ \forall (e,f)\in (\N^*)^2, \ e\wedge f = 1,\ (em,fm\in S \implies m\in S)$.
    %\item $ \forall m\in M,\ \forall k\in \N^*\ (km,(k+1)m\in S \implies m\in S)$.
    \item $ \forall m\in M, \ 2m,3m\in S \implies m\in S$.
    \item $ \forall m\in M,\ (\exists k\in \N, \ \forall k'>k, k'm\in S) \implies m\in S$.
    \end{enumerate}
\end{Lemme}

\begin{comment}
\begin{proof}
    $1\implies 2\implies 3$ are straightforward.
    
    $3 \implies 4 :$ Let $m\in M$. Suppose there exists $k\in \N$ such that, for all $k'>k$, we have $k'm\in S$. Without loss of generality, we can suppose $k$ minimal. If $k>2$, then $2(k-1)m, 3(k-1)m \in S$ so, by 3), we have $(k-1)m \in S$ which is in contradiction with the minimality of $k$. Then $k=2$ and so $m\in S$.
    
    $4 \implies 1 :$ It suffices to prove that if $e,f \in \N$ such that $e \wedge f=1$ then there exists $k \in \N$ such that 
    \[
    \{ n \in \N \mid n \geq k\} \subset \N e+\N f.
    \]
    Since $e$ and $f$ are coprime, we can fix $u,v \in \Z$ such that $ue+vf=1$. Without loss of generality, we can suppose $u\geq 0$ and $v\leq 0$.\\ Let $k\coloneqq -vef$. Then for any $\ell \in \{0,\ldots,e-1\}$, we have
    \[
    k+\ell=-vef+\ell ue+\ell vf=\ell ue-v(e-\ell)f \in \N e+\N f
    \]
    Hence, for any $p=ae+\ell \in \N$, we have
    \[
    k+p=k+\ell+ae \in \N e+\N f
    \]
\end{proof}
\end{comment}

\begin{Def}
    A monoid $S \subset M$ is called "semisaturated" if it verifies one of the properties of \Cref{LemmaSemisaturatedEquivalences}.
\end{Def}

\begin{Ex}\label{ExempleNotSemisaturated}
    The monoid $\Gamma = \N(0,1)+\N(1,2)+\N(2,0)$ (see Figure \ref{fig:monoidNotSemisaturated}) is not semisaturated otherwise $(1,1)$ should belong to it since $(2,2),(3,3)\in \Gamma$.
\begin{figure}[h]
    \centering
    \begin{tikzpicture}[scale=1]
    \tikzmath{ integer \nb; \nb = 4;}
    
    \tikzmath{ integer \nbd; \nbd = \nb-1;}
        \foreach \i in {0,...,\nb}{%
            \foreach \j in {0,...,\nb}{%
            \node at (\i, \j) {$\bullet$};
            }
        }
        \node at (1,1){\color{white} $\bullet$};
        \node at (1,0){\color{white} $\bullet$};
        \node at (3,1){\color{white} $\bullet$};
        \node at (3,0){\color{white} $\bullet$};

        \node at (1,1){$\times$};
        \node at (1,0){$\times$};
        \node at (3,1){$\times$};
        \node at (3,0){$\times$};

    \draw[dashed, -] (0,0) -- (4.25,4.25);
    \draw[->] (0,0) -- (0,0.92);
    \draw[->] (0,0) -- (0.95,1.92);
    \draw[->] (0,0) -- (1.92,0);
    \end{tikzpicture}
    \caption{Illustration of \Cref{ExempleNotSemisaturated}}
    \label{fig:monoidNotSemisaturated}
\end{figure}   
\end{Ex}

\begin{Def}
    Let $S$ be a monoid. The interior of $S$ is defined as 
    \[ \Int(S) \coloneqq \Big\{ x\in S \mid \forall y\in S, \exists n\in \N, \exists z\in S, nx=y+z \Big\} \]
\end{Def}

\begin{Rem}
    If $S \subset M$ is a monoid associated to a cone $\sigma$, then $\Int(S) = S\cap \Int(\sigma)$ where $\Int(\sigma)$ is the relative topological interior of $\sigma$.
\end{Rem}

Let us state the following small lemma, which will be used later.
\begin{Lemme}\label{LemSomInt}
    Let $S$ be a monoid. Then $S + \Int(S) = \Int(S)$.
\end{Lemme}

\begin{proof}
    The inclusion $\Int(S) \subset S + \Int(S)$ is obvious. To prove the converse inclusion, let us consider $x\in \Int(S)$ and $y\in S$. Then, for all $z\in S$, there exists $n\in \N$ such that $nx-z\in S$. Then $n(x+y)-z = (nx-z)+ny \in S$. So $x+y\in \Int(S)$.
\end{proof}

The following property on the interior of monoids clarifies the reason for their link with semisaturation.

\begin{Lemme}[\cite{Swan1992Gubeladze}, Lemma 6.5]\label{LemInterieurNormal}
    Let $S\subset M$ be a monoid and $m\in M$. Then 
    $$\exists k\in \N, km\in \Int(S) \implies k'm\in \Int(S), \text{ for every $k'$ large enough.}$$
\end{Lemme}

\begin{proof}
    Let $m\in M$ such that, there exists $k\in \N$, $km\in \Int(S)\cup\{0\}$. We write $w = km$. We can suppose that $w\in \Int(S)$ because if $w=0$, then $m=0$. Since $m\in M$, we can consider $s_1,s_2\in S$ such that $m=s_1-s_2$. By the definition of $\Int(S)$, we have 
    $$ \forall r\in\{0,\dots,k-1\}, \exists k_r\in\N,\text{\hspace{0.3cm}} k_rw-rs_2\in S $$
    So $k_rw+rm = (k_rw-rs_2)+rs_1 \in S$. Then, for every $n\in \N$, we can do the Euclidean division of $n$ by $k$ and we get $n = qk+r$ with $r\in\{0,\dots,k-1\}$. Then 
    $$ (k_rk+n)m = k_rw+nm = k_rw+(qkm+rm) =(k_rw+rm)+qw \in S $$
    This means that, for all $k'\geqslant k_rk$, we have $k'm\in S$ and, by \cref{LemSomInt}, we get $k'm\in \Int(S)$ if $k'> k_rk$.
\end{proof}

\begin{Prop}
    Let $S\subset M$ be a semisaturated monoid. Then $\Int(S)$ is saturated.
\end{Prop}

\begin{proof}
    This follows by definition from \Cref{LemInterieurNormal}.
\end{proof}

\begin{Thm}[\cite{REID2001703}, Theorem 4.3]\label{TheoReidRoberts}
    Let $S$ be a finitely generated submonoid of $\Z^n$. Then, its semisaturation is
    \[
    S^+:=\bigcup_{\tau \prec \R_{\geq 0} S} \Z(S \cap \tau) \cap \Int(\tau)
    \]
    where $\tau \prec \R_{\geq 0} S$ if $\tau$ is a face of $\R_{\geq 0} S$.
\end{Thm}

Note that, Reid and Roberts proved this theorem for submonoids of $\N^n$, but the proof works for submonoids of $\Z^n$. Moreover, let us remark that the union is disjoint by the following classical lemma
\begin{Lemme}[\cite{fulton} Section 1.2, 7]\label{LemDecompInt}
    Let $\sigma$ be a cone, then $\sigma = \coprod_{\tau\prec\sigma} \Int(\tau)$.
\end{Lemme}

\begin{Ex}\label{ExempleSemisaturated}
    By \Cref{TheoReidRoberts}, the seminormalization of the monoid $\Gamma$ of \Cref{ExempleNotSemisaturated} is obtained by normalizing the interior of each face. Thus $\Gamma^+$ is the monoid generated by $(0,1),(1,1)$ and $(2,0)$, see Figure \ref{fig:monoidSemisaturated}.
    \begin{figure}[h]
    \centering
    \begin{tikzpicture}[scale=1]
    \tikzmath{ integer \nb; \nb = 4;}
    
    \tikzmath{ integer \nbd; \nbd = \nb-1;}
        \foreach \i in {0,...,\nb}{%
            \foreach \j in {0,...,\nb}{%
            \node at (\i, \j) {$\bullet$};
            }
        }
        
        \node at (1,0){\color{white} $\bullet$};
        \node at (3,0){\color{white} $\bullet$};

        \node at (1,0){$\times$};
        \node at (3,0){$\times$};

    \draw[->] (0,0) -- (0,0.92);
    \draw[->] (0,0) -- (0.95,0.92);
    \draw[->] (0,0) -- (1.92,0);
    \end{tikzpicture}
    \caption{Illustration of \Cref{ExempleSemisaturated}}
    \label{fig:monoidSemisaturated}
\end{figure}
\end{Ex}

\begin{Not}
%\todo{Définition à changer}
    Let $\sigma \subset \R^d$ be a cone and $M\subset \R^d$ be a lattice. Let $\{G^\tau\}_{\tau\prec\sigma}$ be a collection of groups such that $G^\sigma = M$ and, for all $\tau\prec\tau'$, we have $G^{\tau}$ is a finite index subgroup of $G^{\tau'}\cap \Span(\tau)$.
    We define 
    \[
    \Gamma^\sigma \coloneqq \coprod_{\tau \prec \sigma} G^\tau\cap \Int(\tau).
    \]
\end{Not}

\begin{Rem}
    The fact that we are considering finite index subgroups in the definition implies that, if $\tau$ is a $d$-dimensional face of $\sigma$, then $G^\tau$ is of rank $d$.
\end{Rem}

We prove now some combinatorial lemmas that we will need later in the paper. In particular, those lemmas prove \cite[Remark 4.4]{REID2001703} which says that the construction defined above always produces a semisaturated monoid.

\begin{Lemme}\label{LemConeTranslate}
    Let $\tau\subset \R^d$ be a cone of dimension $d$ and $G^\tau$ its associated group. Then 
    \[
    \forall g\in G^\tau, \ (g+\tau)\cap\tau\cap G^\tau \neq \emptyset
    \]
\end{Lemme}

\begin{proof}
    Since $\tau$ is of dimension d, we can find a simplicial cone $\Tilde{\tau} = \Cone(e_1,...,e_d) \subset \tau$ such that $\Span_\R(e_1,\dots,e_d) = \R^d$.
    %Since $\tau$ is simplicial, we have $\tau = \Cone(e_1,...,e_d)$ such that $\Span_\R(e_1,\dots,e_d) = \R^d$.
    Let $g\in G^\tau \subset \R^d$, then we can write $g=\sum_i \lambda_ie_i$. We define $g' = \sum_i \max(0,-\lambda_i)e_i \in \tau$. Then $g+g' \in (g+\tau)\cap\tau$.

    Since $G^\tau$ is of rank d, we have $\Q G^\tau = \Span_\Q(e_1,\dots,e_d) \eqqcolon \Tilde{\tau}_\Q$ which is dense in $\Tilde{\tau}$. So there exists an element $\Tilde{g} = \Tilde{\tau}_\Q\cap(g+\tau)\cap\tau \neq \emptyset$. Since $(g+\tau)\cap\tau$ is stable by addition, we can multiply $\Tilde{g}$ by a large enough natural number $n$ until $n\Tilde{g}\in G^\tau \cap (g+\tau)\cap\tau$.
\end{proof}

We first need to show that each $G^\tau$ is generated by its restriction to $\Int(\tau)$.

\begin{Lemme}\label{Lem_GtauGeneParInt}
    For all $\tau\prec\sigma$, we have $G^\tau = \Z(\Int(G^\tau\cap\tau))$.
\end{Lemme}

 \begin{proof}
     The inclusion $\supseteq$ is clear. In order to prove the other inclusion, let us show that, for any $x\in \Int(G^\tau\cap \tau)$, we have :
    \[
    (*)\text{ }\forall g\in G^\tau, \exists n\in \N,\text{ \hspace{0.2cm}} g+nx\in G^\tau\cap \tau.
    \]
    By \Cref{LemConeTranslate}, we can consider $y\in (g+(G^\tau\cap\tau))\cap (G^\tau\cap\tau)$. In particular, we have $g-y\in G^\tau\cap\tau$. So, by definition of $x\in \Int(G^\tau\cap\tau)$, we can consider $n\in \N$ such that $nx-(y-g)\in G^\tau\cap\tau$. Then $nx+g\in G^\tau\cap\tau$ because $y\in G^\tau\cap\tau$ by assumption. We have proven $(*)$.
    Now, let us consider a set of generators $g_1,\dots,g_d$ of $G^\tau$. We apply $(*)$ to every $g_i$ and so we get, for all $i\in\{1\dots,d\}$, that 
    \[
    \exists n_i\in \N,\text{ \hspace{0.2cm}}g_i+n_ix\in G^\tau\cap\tau
    \]
    %Since every $g_i\in G^\tau$ and $x\in G^\tau$, we get $g_i+n_ix\in G^\tau\cap\tau$.
    By \Cref{LemSomInt}, since $x\in \Int(G^\tau\cap\tau)$, we get $g_i+(n_i+1)x\in \Int(G^\tau\cap\tau)$. Then the set $\{x\}\cup \{g_i+(n_i+1)x\}_i \subset \Int(G^\tau\cap\tau)$ generates $G^\tau$. So, we get $\Z(\Int(G^\tau\cap\tau))\subset G^\tau$.
 \end{proof}

We can now deduce the following Lemma, which implies that $\Gamma^\sigma$ satisfies the condition of \Cref{TheoReidRoberts}, and so is semisaturated.

\begin{Lemme}\label{Lem_GammaInterTauGenereGtau}
    For all $\tau\prec\sigma$, we have $G^\tau = \Z(\Gamma^\sigma \cap \tau)$.
\end{Lemme}

\begin{proof}
    Let $\tau\prec\sigma$. By \Cref{LemDecompInt}, we have 
    \[
    \tau = \coprod_{\theta\prec \tau\prec \sigma} \Int(\theta)
    \]
    and so we get 
    \[
    \Gamma^\sigma \cap \tau = \Big(\coprod_{\theta \prec \sigma}G^\theta\cap\Int(\theta) \Big)\cap \tau = \coprod_{\theta \prec \tau} G^\theta \cap \Int(\theta)
    \]
    Since $\theta\prec\tau$ implies $G^\theta\subset G^\tau$, we have $\Gamma^\sigma \cap\tau \subset G^\tau$ and, because $G^\tau$ is a group, we get $\Z(\Gamma^\sigma \cap \tau)\subset G^\tau$.
    Conversely, observe that 
    \[
    G^\tau\cap\Int(\tau)\subset \coprod_{\theta\prec\tau\prec\sigma} G^{\theta}\cap\Int(\theta) = \Gamma^{\theta}\cap\tau
    \]
    So, by \Cref{Lem_GtauGeneParInt}, we get $G^\tau = \Z(G^\tau\cap\Int(\tau))\subseteq \Z(\Gamma^\sigma \cap \tau)$.
\end{proof}

\section{Seminormal affine toric varieties}\label{SectionSeminormalAffTorVar}

In this section, we recall the notion of seminormalization of an algebraic variety and discuss the role of normality in toric geometry. We then explain how seminormality manifests itself combinatorially for affine toric varieties.

\begin{Def}
    Let $A \inj B$ be an integral extension of rings. The seminormalization of $A^+_B$ of $A$ in $B$ is the largest subextension $A\inj C\inj B$ such that:
    \begin{itemize}
        \item $\Spec(C) \to \Spec(A)$ is bijective.
        \item For every $\mathfrak{q} \in \Spec(C)$, then $\kappa(\mathfrak{q}) \to \kappa(\mathfrak{q}\cap A)$ is an isomorphism.
    \end{itemize}
    The ring $A$ is seminormal in $B$ if $A = A^+_B$.
\end{Def}

\begin{Def}
    Let $X$ be an affine variety over a field $k$. Then $X$ is seminormal if its coordinate ring $k[X]$ is seminormal in the coordinate ring $k[X']$ of its normalization. A non-affine algebraic variety is seminormal if it is seminormal on every affine open set.
\end{Def}

If $k$ is an algebraically closed closed field of characteristic $0$, then the seminormalization $X^+$ of $X$ is the largest variety admitting morphisms $X' \to X^+ \to X$ such that $X^+\to X$ is bijective (see \cite{LeahyVitulli} Theorem 2.2). As explained in the introduction, we have the following equivalent definitions of seminormality.

\begin{Prop}[\cite{LeahyVitulli}]\label{PropEquivalenceSeminormalite}
Let $X$ be an affine variety over a field $k$. The following are equivalent:
\begin{enumerate}
    \item $X$ is seminormal
    \item $\forall f\in K(X),\text{ \hspace{0.1cm}}f^2,f^3\in \k[X] \implies f\in \k[X]$.
    \item $\forall f\in K(X),\text{ \hspace{0.1cm}} \Big( \exists N\in \N, \forall n>N, \text{ \hspace{0.1cm}}f^n\in \k[X]\Big) \implies f\in \k[X]$.
\end{enumerate}
\end{Prop}

The seminormalization of a complex algebraic variety can also be obtained, as in the analytic case, by replacing its structural sheaf by the sheaf of rational functions which extend continuously, for the strong topology, on the closed points of the variety, see \cite{Bernard2021} or \cite{BFMQ} for algebraic closed fields of characteristic $0$.

We are interested in studying seminormal toric varieties. Let us start with the affine case. Recall that an affine toric variety is an irreducible algebraic variety with a torus acting on it with a dense orbit. They can also be described as an irreducible variety whose associated ideal is a "toric" ideal i.e. a prime ideal generated by binomials. The main feature of these varieties is their combinatorics description :

Let $V$ be an affine toric variety with torus $T$ acting on it. Note $M$ the (free abelian) group of characters of $T$ i.e. $M=\Hom_{\textbf{Alg gp}}(T,\G_{m,k})$. Then there exists a finitely generated submonoid $S$ of $M$ (\cite[Theorem 1.17]{cox}) such that \[\k[V]=\k[S]\]

Conversely, such a monoid $S$ defines a toric variety $\Spec(\k[S])$ (since $S$ is finitely generated, $\k[S]$ is also finitely generated) with the natural action induced by the $\k$-algebra morphism
\[
\k[S] \to \k[S] \otimes_\k \k[S] \to \k[M] \otimes_\k \k[S]
\]
where the first map is the diagonal morphism defined by $\chi^\gamma \mapsto \chi^\gamma \otimes \chi^\gamma$ and the second map is induced by the inclusion morphism $S \inj M$.

\begin{Ex}\label{EternelExempleAlpha}
Let $\alpha \in \N \setminus \{0\}$. Consider the toric variety $V\coloneqq V(x^2-y^\alpha z) \subset \k^3$. Its torus is 
\[
V \cap (\k^*)^3=\left\{ (x,y,z) \in (\k^*)^3 \mid z=\frac{x^2}{y^\alpha}\right\} \simeq (\k^*)^2
\]
Then 
\[
\k[V]=\k[x,y,z]/(x^2-y^\alpha z) \simeq \k[x,y,x^2y^{-\alpha}] 
\]
Hence, if $S=\N (1,0)+\N (0,1)+\N(2,-\alpha) \subset \Z^2$ then $\k[V]=\k[S]$.\\

\begin{figure}[h]
    \centering
    \begin{tabular}{cc|ccc|cc}
    \begin{tikzpicture}[scale=0.89]
    \tikzmath{ integer \nb; \nb = 3;}
    
    \tikzmath{ integer \nbd; \nbd = \nb-1;}
        \foreach \i in {0,...,3}{%
            \foreach \j in {-\nb,...,2}{%
            \node at (\i, \j) {$\bullet$};
            }
        }

    \draw[->] (0,0) -- (0,0.92);
    \draw[->] (0,0) -- (0.92,0);
    \draw[->] (0,0) -- (1.92,-0.92);

    \end{tikzpicture}&&&

    \begin{tikzpicture}[scale=0.89]
    \tikzmath{ integer \nb; \nb = 3;}
    
    \tikzmath{ integer \nbd; \nbd = \nb-1;}
        \foreach \i in {0,...,3}{%
            \foreach \j in {-\nb,...,2}{%
            \node at (\i, \j) {$\bullet$};
            }
        }
        
        \node at (1,-1){\color{white} $\bullet$};
        \node at (3,-3){\color{white} $\bullet$};

        \node at (1,-1){$+$};
        \node at (3,-3){$+$};

    \draw[->] (0,0) -- (0,0.92);
    \draw[->] (0,0) -- (0.92,0);
    \draw[->] (0,0) -- (1.92,-1.92);
    \end{tikzpicture}&&&

    \begin{tikzpicture}[scale=0.89]
    \tikzmath{ integer \nb; \nb = 3;}
    
    \tikzmath{ integer \nbd; \nbd = \nb-1;}
        \foreach \i in {0,...,3}{%
            \foreach \j in {-\nb,...,2}{%
            \node at (\i, \j) {$\bullet$};
            }
        }
        
        \node at (1,-1){\color{white} $\bullet$};

        \node at (1,-1){$+$};

    \draw[->] (0,0) -- (0,0.92);
    \draw[->] (0,0) -- (0.92,0);
    \draw[->] (0,0) -- (1.92,-2.92);
    \end{tikzpicture}\\
    $\alpha = 1$&\multicolumn{5}{c}{$\alpha = 2$}&$\alpha = 3$\\
    \end{tabular}
    
    \caption{Illustration of \Cref{EternelExempleAlpha}}
    \label{fig:EternelExempleAlpha}
\end{figure}
\end{Ex}

In the classical literature, a toric variety is almost always assumed normal in order to get a combinatorial description, which is given by the following Theorem.
\begin{Thm}[\cite{cox} Theorem 1.3.5]\label{TheoDualiteConeVarNormal}

Let $V$ be an affine toric variety, and let $M$ be the group of characters of its torus. The following are equivalent:
\begin{itemize}
    \item $V$ is normal.
    \item $V=\Spec(\k[S])$, $S \subset M$ is a saturated finitely generated monoid.
    \item $V=\Spec(\k[S_\sigma])$  where $S_\sigma=\sigma^\vee \cap M$ and $\sigma$ is a strongly convex rational polyhedral cone.
\end{itemize}
\end{Thm}

Hence, if the variety is normal, we recover the monoid by intersecting the dual cone with the lattice. This means, in a sense, that, for a given lattice, we can forget the information of the monoid and only care about the cones.

\begin{Ex}\label{ExempleDualiteConeVarNormal}
In \Cref{EternelExempleAlpha}, the variety is normal if and only if $\alpha=1$. We recover the monoid $S$ by intersecting the lattice $M$ with the dual of the cone $\sigma = (1,2)\R_{\geq 0} + (1,0)\R_{\geq 0}$. See Figure \ref{fig:NormalConeIntersect}.

\begin{figure}[h]
    \centering
    \begin{tabular}{ccc}
    \begin{tikzpicture}[scale=0.8]
    
    \tikzmath{ integer \nb; \nb = 3;}
    
    \tikzmath{ integer \nbd; \nbd = \nb-1;}
        \foreach \i in {0,...,4}{%
            \foreach \j in {-2,...,3}{%
            \node at (\i, \j) {};
            }
        }

    \draw[-] (0,-1) -- (1.6,2.2);
    \draw[-] (0,-1) -- (4,-1);
    \end{tikzpicture}

    &
    
    \begin{tikzpicture}[scale=0.8]
    \tikzmath{ integer \nb; \nb = 3;}
    
    \tikzmath{ integer \nbd; \nbd = \nb-1;}
        \foreach \j in {-2,...,3}{%
        \node at (1, \j) {};
        }

    \draw[->] (0,0) to [bend left=20] (3,0);
    \end{tikzpicture}
    
    &

    \begin{tikzpicture}[scale=0.8]
    \tikzmath{ integer \nb; \nb = 3;}
    
    \tikzmath{ integer \nbd; \nbd = \nb-1;}
        \foreach \i in {0,...,4}{%
            \foreach \j in {-2,...,3}{%
            \node at (\i, \j) {$\bullet$};
            }
        }

    \draw[->] (0,0) -- (0,0.92);
    \draw[dashed, -] (0,0) -- (0,3.1);
    \draw[->] (0,0) -- (2,-1);
    \draw[dashed, -] (0,0) -- (4,-2);

    %\node at (1.5, 1.5) {$\sigma^{\vee}$};
    \end{tikzpicture}\\
    $\sigma \subset \R^2$ & & $\sigma^{\vee} \cap M$
    \end{tabular}
    
    \caption{Illustration of \Cref{ExempleDualiteConeVarNormal} }
    \label{fig:NormalConeIntersect}
\end{figure}
\end{Ex}

More precisely, we have an equivalence of categories between the category of strongly convex cones with linear morphisms between them and the category of affine toric varieties with equivariant regular morphisms between them.

Let us state now the link between seminormal affine toric varieties and seminormal monoids.

\begin{Thm}[\cite{REID2001703} Theorem 3.8]\label{TheoVarSNiffMonoidSemiSat}
%\todo{Déjà fait dans RR ?}
    Let $V$ be an affine toric variety and $M$ be its character lattice. Then $V$ is seminormal if and only if $V=\Spec(\k[S])$ where $S\subset M$ is semisaturated.
\end{Thm}

\begin{proof}
    Let us suppose that $V$ is a seminormal toric variety with coordinate ring $\k[S]$. Let $m\in M$ be such that $2m\in S$ and $3m\in S$. We have $\chi^m\in \k[M] \subset K(V)$ and $(\chi^m)^2, (\chi^m)^3\in \k[S]$, so, by \Cref{LemmaSemisaturatedEquivalences}, we obtain $\chi^m\in \k[S]$ i.e. $m \in S$.

    If $S$ is semisaturated then by \cite[Theorem 3.8]{REID2001703}, we have
    \[
    \k[V^+]\coloneqq \k[V]^+_{\k[V]'}=\k[S]^+_{\k[S]'}=\k[S^+]
    \]
    Since $S$ is semisaturated then $S=S^+$ and
    \[
    \k[V^+]=\k[S]=\k[V]
    \]
    i.e. $V$ is seminormal.
    %{\color{blue} C'est vrai cette réciproque ?} Conversely, let $f\in K(V)$ such that $f^2,f^3\in \k[S]$. {\color{red}\textsc{bullshit}}Then there exists $m\in M$ such that $f = \chi^m$ and so $\chi^{2m}, \chi^{3m}\in \k[S]$. Then $2m,3m\in S$ implies that $m\in S$, since $S$ is saturated, and so $f = \chi^m \in \k[S]$. 
\end{proof}%

\begin{Ex}
    The case $\alpha = 3$ of \Cref{EternelExempleAlpha}, given by the equation $x^2-y^{3}z$, is not seminormal because $x/y$ is not polynomial while $(\frac{x}{y})^2 = yz$ and $(\frac{x}{y})^3 = xz$ are polynomial. We can see in Figure \ref{fig:EternelExempleAlpha} that its associated monoid is not semisaturated because $(1,-1)$ belongs to the interior of the saturated maximal face. However, $x^2-y^{2}z$ is seminormal because the interior of each face is equal to the interior of the saturated face. To conclude, the variety $x^2-y^{\alpha}z$ is normal if and only if $\alpha \leq 1$ and is seminormal if and only if $\alpha \leq 2$.
\end{Ex}

\section{General seminormal toric varieties}\label{SectionGeneralTorVar}

In this section, we prove that seminormal toric varieties correspond to a combinatorial object that we call "fan with attached groups". In order to do so, we use the theory of non-normal toric varieties introduced by Teissier and González Peréz in \cite{TeissierGonzales2012}.

%\subsection{Construction of González Peréz and Teissier}
\subsection{Non-normal toric varieties and fans with attached monoids}\label{SubsecTeissierGonzPerez}

Let us recall briefly, in this subsection, the construction of \cite{TeissierGonzales2012}.

\begin{Def} \label{def_fan-GPT}
A fan with attached monoids\footnote{In the original article, it is called "fan with attached semigroups".} is a triple $(N, \Sigma,
\Gamma)$ consisting of lattice $N$,  a fan $\Sigma$ in $N_{\R}$ and a family of finitely generated submonoids $\Gamma = \{
\Gamma_\sigma \subset \sigma^\vee \cap M\}_{\sigma \in \Sigma}$ of a lattice $M =
\Hom (N, \Z)$ such that:
\begin{enumerate}
\item  $\Z \Gamma_\sigma = M$ and $\R_{\geq 0}\Gamma_\sigma=\sigma^\vee$, for $\sigma  \in \Sigma$.
\item $\Gamma_\tau = \Gamma_\sigma + M (\tau, \Gamma_\sigma)$, for
each $\sigma \in \Sigma$ and  any face $\tau$ of $\sigma$.
\end{enumerate}
where $M (\tau, \Gamma_\sigma)$ is the free abelian group generated by $\Gamma_\sigma \cap \tau^\perp$.
\end{Def}

 The  corresponding toric variety ${T}_{\Sigma}^{\Gamma}$
 is the union of the affine varieties ${T}^{\Gamma_\sigma} \coloneqq \Spec(\k[\Gamma_\sigma])$
for $\sigma \in \Sigma$ where for any pair $\sigma, \sigma'$ in $\Sigma$ we
glue up ${T}^{\Gamma_\sigma}$ and ${T}^{\Gamma_{\sigma'}}$ along their common open affine variety ${T}^{\Gamma_{\sigma \cap \sigma'}}$ (the inclusion $\Gamma_\sigma \subset \Gamma_{\sigma \cap \sigma'}$ induces an open immersion ${T}^{\Gamma_{\sigma \cap \sigma'}} \inj {T}^{\Gamma_{\sigma}}$ cf. \cite[Lemma 3.10]{TeissierGonzales2012}).

\begin{Ex}
   Let $\sigma_1 = (\R_{\geq 0})^2$ with $\Gamma_{\sigma_1} = \N \{(2,0), (3,0),(1,1), (0,2)\}$ and $\sigma_2 = \R_{\geq 0}\times\R_{\leq 0}$ with $\Gamma_{\sigma_2} = \N \{ (1,0), (1,-1),(0,-2)\}$. Let $\tau = \sigma_1\cap\sigma_2$. Then $M(\tau, \Gamma_{\sigma_1})= M(\tau, \Gamma_{\sigma_2})$ and yet condition ii) of \Cref{def_fan-GPT} is not verified.
\end{Ex}

\begin{Def}[\cite{TeissierGonzales2012}, Definition 6.1]
   A morphism of fans with attached monoids $(N,\Sigma,\Gamma)\to (N',\Sigma',\Gamma')$ is a morphism $\varphi : N \to N'$ such that for every $\sigma \in \Sigma$, there exists $\sigma' \in \Sigma'$ such that $\varphi^\top(\Gamma_{\sigma'}) \subset \Gamma_{\sigma}$.
   We denote by $\mathfrak{Fan}^{mon}$ the category of fans with attached groups.
\end{Def}

\begin{Def}[\cite{TeissierGonzales2012}, Definition 7.4]
    An action of a group on an algebraic variety is \textit{good} if $X$ is covered by a finite number of affine open subsets which are invariant by the action.
\end{Def}

\begin{Def}[\cite{TeissierGonzales2012}, Definition 7.5]
    A toric variety $X$ is an irreducible separated algebraic variety equipped with a good action of an algebraic torus $T$ embedded in $X$ as a Zariski open set such that the action of $T$ on $X$ extends the action of $T$ on itself by multiplication.
    We note $\mathfrak{TorVar}$ the category of such varieties with equivariant regular morphisms.
\end{Def}

\begin{Rem}
    As shown in \cite[Example 7.3]{TeissierGonzales2012}  with the nodal cubic $C\subset \P^2_\k$ given by the equation $y^2z-x^2(x+z)$, the assumption of having a "good" action cannot be dropped if we remove the normality (see \cite{Sumihiro}) assumption in the definition of a toric variety. Since $C$ is a curve with only normal crossing singularities, it is seminormal. So, this example also shows that we cannot drop the assumption of having a "good" action for seminormal toric varieties.
\end{Rem}

\begin{Thm}[\cite{TeissierGonzales2012} Corollary 7.7] \label{EquivalenceCat_GPT}
    The functor $(N,\Sigma,\Gamma) \mapsto T_\Sigma^\Gamma$ is an equivalence of categories between the category of fans with attached monoids $\mathfrak{Fan}^{mon}$ and the category of toric varieties with a good action $\mathfrak{VarTor}$.
\end{Thm}

\subsection{Fans with attached groups and seminormal toric varieties}\label{SubsecFansWithAttachedGroups}

We now present a construction similar to the one of González Pérez and Teissier, but where the semisaturation property of the monoids allow for simpler gluing conditions. 

\begin{Def}\label{DefFanAttachedGroups}
    A fan with attached groups is given by the datum of a triple $(N,\Sigma, G)$ consisting of a lattice $N$, a fan $\Sigma$ in $N_{\R}$ and a collection of subgroups $G = \{G_\sigma\}_{\sigma\in\Sigma}$ of a lattice $M=\Hom(N,\Z)$ such that
    \begin{itemize}
        \item[1)] $G_0 = M$
        \item[2)] If $\tau\prec\sigma$, then $G_\sigma$ is a finite index subgroup of $G_\tau \cap \sigma^\perp$.
    \end{itemize}
    %\todo{Je pense que l'on peut remplacer $M \cap \sigma^\perp$ par $\sigma^\perp$ }
\end{Def}

\begin{Ex} %\todo{Je crois qu'il y a un problème dans l'exemple}
    Let $\Sigma$ be the fan represented in Figure \ref{fig:ExampleFanAttachedGroups}. We show how to associate a collection of groups to $(\Z^2, \Sigma)$. We start with the cones of lowest dimension. For the $0$-dimensional cone $0$, we have no choice: $G_0 = \Z^2$. For the one dimensional cones, the condition is $G_{\tau_i} \subseteq \Z^2\cap \tau_i^{\perp}$. So, for example, we can choose
    $$\begin{cases}
        G_{\tau_1} = (1,2).\Z\\
        G_{\tau_{1,2}} = (0,1).2\Z\\
        G_{\tau_{2,3}} = (2,1).2\Z\\
        G_{\tau_3} = (0,1).3\Z
    \end{cases}.$$
    The groups $G_{\sigma_i}$ must all be $\{0\}$ because, since they are of maximal dimension, we have $\sigma_i^{\perp} = \{0\}$.

    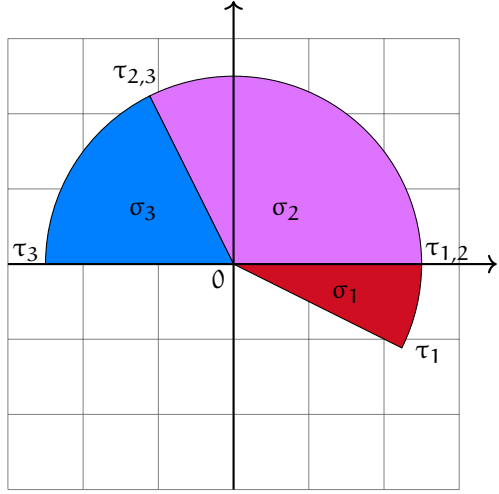
\begin{figure}[h]
    \centering
    \begin{tikzpicture}[scale=1]

    % Draw grid
    \draw[step=1cm,gray,very thin] (-3,-3) grid (3,3);
    
    \tikzmath{\ray = 2.5;}
    \filldraw[fill opacity=0.5,fill=couleur1] (0,0) -- (\ray,0) arc (0:116.5:\ray cm) -- cycle;
    \filldraw[fill opacity=0.5,fill=couleur3] (0,0) -- (-\ray,0) arc (180:116.5:\ray cm) -- cycle;
    \filldraw[fill opacity=0.5,fill=couleur2] (0,0) -- (\ray,0) arc (360:333.5:\ray cm) -- cycle;

    \node at (1.5,-0.4) {$\sigma_1$};
    \node at (0.7,0.7)  {$\sigma_2$};
    \node at (-1.2,0.7) {$\sigma_3$};

    \node at (-0.2,-0.2) {$0$};

    \node at (2.6,-1.2) {$\tau_{1}$};
    \node at (2.85,0.15) {$\tau_{1,2}$};
    \node at (-1.3,2.5) {$\tau_{2,3}$};
    \node at (-2.75,0.15) {$\tau_{3}$};

    \draw[thick,->] (-3,0) -- (3.5,0) node[anchor=north west] {};
    \draw[thick,->] (0,-3) -- (0,3.5) node[anchor=south east] {};
    
    \end{tikzpicture}
    \caption{A fan with attached groups}
    \label{fig:ExampleFanAttachedGroups}
    \end{figure}
\end{Ex}

\begin{Def}
   A morphism of fans with attached groups $(N,\Sigma,G)\to (N',\Sigma',G')$ is a group morphism $\varphi : N \to N'$ satisfying the following properties :
   \begin{itemize}
       \item[1)] For every $\sigma \in \Sigma$, there exists $\sigma' \in \Sigma'$ such that $\varphi^\top(G_{\sigma'}) \subset G_{\sigma}$.
       \item[2)] For every maximal cone $\sigma$, there exists $\sigma'$ in $\Sigma'$ such that $\varphi(\sigma) \subset \sigma'$.
   \end{itemize}
   We denote by $\mathfrak{Fan}^{gr}$ the category of fans with attached groups.
\end{Def}

\begin{Rem}
    Contrary to the definition of morphisms of fans with attached monoid given in \cite{TeissierGonzales2012}, we must require that every maximal cone of $\Sigma$ is sent to a cone of $\Sigma'$. Let us consider, for example, the fan with attached groups $(\Z^2,\Sigma_{\R_{\geq 0}^2},G)$ define by $\Sigma_{\R_{\geq 0}^2} = \big\{0,0\times \R_{\geq 0},\R_{\geq 0}\times 0,\R_{\geq 0}^2\big\}$ and $G = \big\{\Z^2,\Z\times 0,0\times \Z, 0\big\}$. Then we consider the morphism of groups $\varphi : \Z^2 \to \Z^2, (x,y)\mapsto (-x,-y)$. For every $\sigma\in \Sigma_{\R_{\geq 0}^2}$, we have $\varphi^\top(G_{\sigma})=G_{\sigma}$. However, it does not send the cone $\R_{\geq 0}^2$ in $\R_{\geq 0}^2$, so $\varphi$ is not a morphism of fans. On the contrary, a morphism of fans with attached groups is a fan morphism between the underlying fans.
\end{Rem}

We know from \Cref{EquivalenceCat_GPT} that there is an equivalence of categories between $\mathfrak{Fan}^{mon}$ and $\mathfrak{TorVar}$. By restriction and by \Cref{TheoVarSNiffMonoidSemiSat}, we also have an equivalence of categories between the category $\mathfrak{Fan}^{mon,sn}$ of fans with attached semisaturated monoids and the category $\mathfrak{TorVar}^{sn}$ of seminormal toric varieties with a good action.

\[\begin{tikzcd}
	{\mathfrak{Fan}^{mon}} & {\mathfrak{TorVar}} \\
	{\mathfrak{Fan}^{mon,sn}} & {\mathfrak{TorVar}^{sn}}
	\arrow["\simeq", from=1-1, to=1-2]
	\arrow[hook, from=2-1, to=1-1]
	\arrow[hook, from=2-2, to=1-2]
	\arrow["\simeq", from=2-1, to=2-2]
\end{tikzcd}\]

\begin{Def}
    We define the functor $F : \mathfrak{Fan}^{gr} \to \mathfrak{Fan}^{mon,sn}$ such that, for all fan with attached groups $(N,\Sigma,G)$, then $F(N,\Sigma,G)=(N,\Sigma,\Gamma)$ where $\Gamma = \{\Gamma_{\sigma} \}_{\sigma\in\Sigma}$ with
    \[
    \Gamma_\sigma \coloneqq \coprod_{\tau \preceq \sigma} G_\tau \cap \Int(\sigma^\vee \cap \tau^\perp)
    \]
    and such that, for every morphism $\varphi : N \to N'$, then $F(\varphi)=\varphi : N \to N'$.
\end{Def}

\begin{Not}
    In what follows, if $\tau$ is a face of a cone $\sigma$,  we write $\tau^*$ the cone $\sigma^\vee \cap \tau^\perp$. Recall that the map $\tau \mapsto \tau^*$ is a inclusion-reversing bijection between the faces of $\sigma$ and the faces of $\sigma^\vee$ (cf. \cite[Section 1.2 10)]{fulton})
\end{Not}

\begin{Prop}\label{PropWellDefFunc}
    The functor $F$ is well-defined.
\end{Prop}
 
In the following, we will consider a fan with attached groups $(N,\Sigma,G)$ and the corresponding family $\Gamma = \{\Gamma_{\sigma} \}_{\sigma\in\Sigma}$ such that $F(N,\Sigma,G)=(N,\Sigma,\Gamma)$.

As a first step to prove that $F$ is well-defined, we state the following proposition, which is a direct consequence of \Cref{Lem_GammaInterTauGenereGtau}.

\begin{Prop}\label{Prop_LesGCorrespondentAuxM}
    Let $\tau\prec\sigma$ be cones of $\Sigma$, then 
    \[
    M(\tau,\Gamma_\sigma) = G_{\tau}
    \]
\end{Prop}

\begin{proof}
    Since $G^{\tau^*} = G_\tau$, we simply apply \Cref{Lem_GammaInterTauGenereGtau} to $\Gamma_{\sigma}$ and $\tau^*$. We get $G^{\tau^*} = \Z(\Gamma_{\sigma}\cap(\sigma^\vee \cap \tau^\perp)) = \Z(\Gamma_{\sigma}\cap \tau^\perp) = M(\tau,\Gamma_\sigma)$.
\end{proof}

We give two technical lemmas that will be needed.

\begin{Lemme}\label{Lem_GammaPlusGSemisaturated}
    Let $\sigma$ be a cone of dimension $d+1$, let $\tau$ be a face of $\sigma$ of codimension $1$ and assume that $\Gamma_{\sigma}$ is semisaturated. Then $\Gamma_\sigma +G_\tau$ is semisaturated.
\end{Lemme}

\begin{proof}
    The group $G_\tau=\Z(\Gamma_\sigma \cap \tau^\perp)$ has rank $1$, meaning it is of the form $\Z v$ for some $v \in \Gamma_\sigma \cap \tau^\perp$ because it is a discrete subgroup of a one dimensional $\R$-vector space. Hence
    \[
    \k[\Gamma_\sigma+G_\tau]=\k[\Gamma_\sigma][v^{-1}].
    \]
    Since $\Gamma_\sigma$ is semisaturated, then $\k[\Gamma_\sigma]$ is seminormal (cf. \Cref{TheoVarSNiffMonoidSemiSat}) and so $\k[\Gamma_\sigma+G_\tau] = \k[\Gamma_\sigma][v^{-1}]$ is seminormal by \cite[Corollary 2.2]{GrecoTraverso}. In conclusion, $\Gamma_\sigma +G_\tau$ is semisaturated.
\end{proof}

\begin{Lemme} \label{Lem_interieur_sur_interieur}
    Let $\Sigma$ be a fan of $N_\R$, $\Sigma'$ be a fan of $N_\R'$ and $\varphi \colon N \to N'$ be a fan morphism. If $\tau \in \Sigma$ and $\tau'$ is the minimal cone of $\Sigma'$ containing $\varphi(\tau)$ then $\varphi(\Int(\tau)) \subset \Int(\tau')$.
\end{Lemme}

\begin{proof}
    Let $\tau \in \Sigma$ and $\tau' \in \Sigma'$ such that $\tau'$ is the minimal cone of $\Sigma'$ containing $\varphi(\tau)$.\\
    Since $\varphi(\mathrm{span}(\tau))=\mathrm{span}(\varphi(\tau))$ (by linearity) then by the open mapping theorem, the restriction of $\varphi_{\mid \mathrm{span}(\tau)} : \mathrm{span}(\tau) \to \mathrm{span}(\varphi(\tau)) $ is an open mapping. Hence,  $\varphi(\Int(\tau)) \subset \Int(\varphi(\tau))$. Thanks to 
    \cite[Lemma 5.2]{Swan1992Gubeladze}
    and the minimality of $\tau'$, we have 
    \[
    \varphi(\Int(\tau)) \subset \Int(\tau').
    \]
\end{proof}

We now prove  \Cref{PropWellDefFunc}.
\begin{proof}
    We have to show that $(N,\Sigma,\Gamma)$ is an element of $\mathfrak{Fan}^{mon, sn}$, i.e. that it respects the conditions of \Cref{def_fan-GPT}. For every $\sigma\in\Sigma$, the set $\Gamma_{\sigma}$ is a semisaturated monoid as we saw in \Cref{SectionSemisaturatedMonoids}.
    By \Cref{Prop_LesGCorrespondentAuxM}, we have
    $$ \Z\Gamma_{\sigma} = \Z(\Gamma_{\sigma}\cap0^\perp) = M(0,\Gamma_{\sigma}) = G_0 = M $$
    So $\Z\Gamma_{\sigma} = M$, for every $\sigma\in\Sigma$. Moreover, $\R_{\geq 0}\Gamma_{\sigma} = \sigma^\vee$ because $\Gamma_{\sigma}$ is non-zero when intersected with any $1$-dimensional face of $\sigma^\vee$. 
    We now prove that, for any face $\tau$ of any cone $\sigma\in \Sigma$, we have
    $$ \Gamma_{\tau} = \Gamma_{\sigma} + G_{\tau}$$
    We start with the case where $\tau$ is a face of codimension $1$. Since $\Gamma_{\tau}$ and $\Gamma_{\sigma} + G_{\tau}$ are semisaturated (cf. \Cref{Lem_GammaPlusGSemisaturated}) and since $\Cone(\Gamma_\sigma+G_\tau)=\Cone(\Gamma_\tau)=\tau^\vee$, it suffices to prove (thanks to \Cref{TheoReidRoberts}), for all $\theta \preceq \tau$ (or, equivalently, for all $\theta^* \preceq \tau^*$), that
    \[
    (\Gamma_\sigma+G_\tau)\cap \Int(\theta^\perp \cap \tau^\vee)=\Gamma_\tau\cap \Int(\theta^\perp \cap \tau^\vee)=G_\theta\cap \Int(\theta^\perp \cap \tau^\vee)
    \]
    Let $\theta$ be a face of $\tau$. Let $x \in \Gamma_\sigma$ and $y \in G_\tau$ such that $x+y \in \Int(\theta^\perp \cap \tau^\vee)$. In particular, we have $x+y\in \theta^\perp$ and $y\in \tau^\perp \subset \theta^\perp$. Then, for all $v\in \theta$, we get
    $\langle x,v\rangle = \langle x+y,v\rangle-\langle y,v\rangle = 0$, and so $x\in \theta^\perp$.
    %This means that $x_j+y_j=0$, for all $j\in J$, and that $y_i=0$, for all $i<d+1$. Since $\theta$ is a face of $\tau$, then $d+1\notin J$ and so we get $x_j=0$, for all $j\in J$. This exactly means that $x\in \theta^\perp$.
    Then $x\in \Gamma_\sigma \cap \theta^\perp \subset  \Z(\Gamma_\sigma \cap \theta^\perp)= G_\theta$. Moreover, we also have $y\in G_\theta$ because $y\in G_{\tau}$ and $\theta \preceq \tau$ implies $G_\tau \subset G_\theta$.
    Hence, $x+y \in G_\theta$ and 
    \[
    (\Gamma_\sigma+G_\tau)\cap \Int(\theta^\perp \cap \tau^\vee)\subset G_\theta\cap \Int(\theta^\perp \cap \tau^\vee)
    \]
    
    Conversely, let $x \in G_\theta \cap \Int(\theta^\perp \cap \tau^\vee)$. Since $\tau$ is a codimension $1$ face of $\sigma$, we can write $\tau = \sigma \cap m^\perp$ where $m$ can be chosen in $G_{\tau}$. Moreover, by \cite[Exercise 1.2.2]{cox}, we have $x\in \Int(\theta^\perp\cap\tau^\vee)$ if and only if $\theta = \tau \cap x^\perp$. So $\theta = \sigma \cap m^\perp \cap x^\perp = \sigma \cap (am+x)^\perp$ for $a\in \N$ large enough in order to get $am+x\in \sigma^\vee$. Hence, we get $am+x \in \Int(\theta^\perp \cap \sigma^\vee)$.
    %Voir démo Fulton (4) ou Oda démo prop A.5.
    Moreover, since $m\in G_\tau \subset G_\theta$, we have $am+x \in G_\theta \cap \Int(\sigma^\vee\cap\theta^\perp)$. Then $am+x \in \Gamma_\sigma \cap \Int(\sigma^\vee \cap \theta^\perp)$ by definition of $\Gamma_\sigma$ and $-am \in G_\tau$. Finally $x = (am+x) - am \in (\Gamma_\sigma+G_\tau)\cap \Int(\theta^\perp \cap \tau^\vee)$.
  
    %Conversely, let $x=(x_1,\ldots,x_{d+1}) \in G_\theta \cap \Int(\tau^\vee \cap \theta^\perp)$. Note $\underline{v}=v e_{d+1}$ a positive generator of $G_\tau$.
    
    %If $x_{d+1}>0$ then $x \in \Int(\sigma^\vee \cap \theta^\perp)$. By definition of $\Gamma_\sigma$, we have $G_\theta \cap \Int(\sigma^\vee\cap\theta^\perp) \subset \Gamma_{\sigma}$. So $x\in \Gamma_\sigma \cap \Int(\sigma^\vee \cap \theta^\perp) \subset (\Gamma_\sigma+G_{\tau}) \cap \Int(\tau^\vee \cap \theta^\perp)$. Now suppose $x_{d+1} \leq 0$. One can consider $\ell \in \N$ such that $x_{d+1}+\ell v>0$. By the previous arguments, we get $x+\ell\underline{v} \in \Gamma_{\sigma}$. Then the decomposition $x=(x+\ell \underline{v})-\ell \underline{v}$ gives us $x \in \Gamma_\sigma+G_\tau$.

    Now, for any $\tau\prec\sigma \in \Sigma$, we can consider $\tau\prec\tau_1\prec\dots\prec\tau_n\prec\sigma \in \Sigma$ where $\tau_i$ is a face of $\tau_{i+1}$ of codimension 1. Since $\tau_i\prec\tau_{i+1}$ implies $G_{\tau_{i+1}}\subset G_{\tau_{i}}$, then we have, for all $i$,
    $$\Gamma_{\tau_i} = \Gamma_{\tau_{i+1}} + G_{\tau_i} = \Gamma_{\tau_{i+2}} + G_{\tau_{i+1}} + G_{\tau_i} = \Gamma_{\tau_{i+2}} + G_{\tau_i}$$
    Hence $$ \Gamma_{\tau} = \Gamma_{\sigma} + G_{\tau}.$$
    Finally, it remains to show that every morphism of fans with attached groups $\varphi : (N,\Sigma,G)\to (N',\Sigma',G')$, is sent, by $F$, to a morphism of fans with attached monoids $(N,\Sigma,\Gamma)\to (N',\Sigma',\Gamma')$. Let $\varphi$ be a morphism of fans with attached groups.\\
    Let $\sigma \in \Sigma$. Because $\varphi$ is supposed to be a morphism of fans, then we can fix $\sigma' \in \Sigma'$ such that $\varphi(\sigma) \subset \sigma'$ (and $\varphi^\top(\sigma'^\vee) \subset \sigma^\vee$). Let $\theta'$ be a face of $\sigma'$. Then $\varphi^\top(\theta'^*) \subset \sigma^\vee$. Note $\tau^*$ the minimal face of $\sigma^\vee$ such that $\varphi^\top(\theta'^*) \subset \tau^*$. Then, by \Cref{Lem_interieur_sur_interieur}, if $x \in \Int(\theta'^*)$ then $\varphi^\top(x)\in \Int(\tau^*)$. Since $\varphi$ is a morphism of fan with attached groups, we can fix a minimal $\tau'\in \Sigma'$ such that $\varphi^\top(G_{\tau'}) \subset G_\tau$. Hence, $\varphi^\top({\tau'}^\perp) \subset \tau^\perp$ and by hypothesis on $\sigma'$, we get $\varphi^\top({\tau'}^*) \subset \tau^*$.
    By minimality, $\tau'$ is a face of $\theta'$ ($\tau'$ is a face of $\sigma'$). We have an inclusion $G_{\theta'} \subset G_{\tau'}$ and hence, $\varphi^\top(G_{\theta'}) \subset G_\tau$.\\
    In conclusion, the map $\varphi^\top$ sends $G_{\theta'} \cap \Int(\theta'^*)$ on $G_{\tau} \cap \Int(\tau^*)$. We can do this construction for any face of $\sigma'$ hence 
    \[
    \varphi^\top\left(\coprod_{\theta' \preceq \sigma'} G_{\theta'} \cap \Int(\theta'^\perp \cap \sigma'^\vee)\right) \subset \coprod_{\tau \preceq \sigma} G_{\tau} \cap \Int(\tau^\perp \cap \sigma^\vee) 
    \]
    i.e. $\varphi^\top(\Gamma_{\sigma'}) \subset \Gamma_\sigma$. The morphism $\varphi$ is a morphism of fans with attached monoids.
\end{proof}

Now that we have proven that the functor $F$ is well-defined, it remains to show that it is an equivalence of categories.

\begin{Thm}[main theorem]
    There is an equivalence of categories between the category $\mathfrak{TorVar}^{sn}$ of seminormal toric varieties with a good action and the category $\mathfrak{Fan}^{gr}$ of fans with attached groups.
\end{Thm}

\begin{proof}
    We start by proving that $F$ is essentially surjective. Let $T^\Gamma_{\Sigma}$ be an element of $\mathfrak{TorVar}^{sn}$ and $(N,\Sigma,\Gamma)$ be its corresponding element in $\mathfrak{Fan}^{mon,sn}$. We consider $(N,\Sigma,\{G_{\sigma}\}_{\sigma\in\Sigma})$ with $G_{\sigma} = M(\sigma,\Gamma_{\sigma})$ for all $\sigma \in \Sigma$. Let $\sigma\in\Sigma$. By assumption, we have that $\Gamma_{\sigma}$ is seminormal and so, by \Cref{TheoReidRoberts},
    \[
    \Gamma_{\sigma} = \coprod_{\tau \preceq \sigma^\vee} \Z(\tau\cap\Gamma_{\sigma})\cap\Int(\tau)
    = \coprod_{\tau \preceq \sigma} \Z(\tau^*\cap\Gamma_{\sigma})\cap\Int(\tau^*).
    \]
    We have $\Z(\tau^*\cap\Gamma_{\sigma}) = \Z(\tau^\perp\cap\Gamma_{\sigma}) = M(\tau,\Gamma_{\sigma})$ and, by \cite[Lemma 3.9, iii)]{TeissierGonzales2012}, we have $M(\tau,\Gamma_{\sigma}) = M(\tau,\Gamma_{\tau})$. Then 
    \[
    \Gamma_{\sigma} = \coprod_{\tau \preceq \sigma} M(\tau,\Gamma_{\tau})\cap\Int(\tau^*)
    = \coprod_{\tau \preceq \sigma} G_{\tau}\cap\Int(\tau^*)
    \]
    We have obtained $F(N,\Sigma,\{G_{\sigma}\}_{\sigma\in\Sigma}) = (N,\Sigma, \Gamma)$ and so $F$ is essentially surjective.
    
    %Pleinement fidèle : Vérifier Hom(A,B) = Hom(F(A),F(B))\todo{C'est la seconde partie de "connaissez-vous les morphismes". ça utilise les Lemmes "Exaircisses".}

    In order to show that $F$ is fully faithful, let us consider a morphism of fans $\varphi$ with attached monoids. Let $\sigma \in \Sigma$. By definition of a morphism of fans with attached monoids, there exists $\sigma'\in \Sigma'$ such that 
    \[
    \varphi^\top(\Gamma_{\sigma'}) \subset \Gamma_\sigma
    \]
    Since $\sigma'^\vee=\R_{\geq 0}\Gamma_{\sigma'}$ (and $\sigma^\vee=\R_{\geq 0}\Gamma_\sigma$) then $\varphi^\top((\sigma')^\vee)\subset \sigma^\vee$.
    This induces two more inclusions by duality : $\varphi(\sigma)\subset \sigma'$ and $\varphi^\top(\sigma'^\perp)\subset \sigma^\perp$. In conclusion, by \Cref{Prop_LesGCorrespondentAuxM}, we get 
    \[
    \varphi^\top(G_{\sigma'})=\varphi^\top(\Z\cdot (\Gamma_{\sigma'} \cap (\sigma')^\perp)) \subset \Z\cdot\varphi^\top(\Gamma_{\sigma'} \cap (\sigma')^\perp) \subset \Z \cdot(\Gamma_\sigma\cap\sigma^\perp)=G_\sigma 
    \]
\end{proof}

\begin{Ex}\label{CollageMonoidesSN}
    Let us consider again the fan $\Sigma$ of Figure \ref{fig:ExampleFanAttachedGroups}. The cone $\sigma_1$, for example, is of maximal dimension, so the attached monoid is 
    \begin{align*}
         \displaystyle \Gamma_{\sigma_1} & = \coprod_{\tau \preceq \sigma_1} G_{\tau}\cap\Int(\tau^*)\\
         & = \Big( G_{0}\cap\Int(\sigma_1^{\vee})\Big) \cup \Big(G_{\tau_{1}}\cap\Int(\tau_{1}^*) \Big)\medskip \cup \Big(G_{\tau_{1,2}}\cap\Int(\tau_{1,2}^*) \Big)\cup \{0\}\\
         & = \Big( \Z^2 \cap\Int(\sigma_1^{\vee})\Big) \cup \Big(\Z\cap (1,2)\N_{>0} \Big)\cup \Big(2\Z\cap(0,-1)\N_{>0} \Big)\cup \{0\}\\
         & = (0,-2)\N + (1,2)\N + (1,1)\N
    \end{align*}
    In the same way, we get that $\Gamma_{\sigma_2}$ is generated by $(2,1);(0,2);(1,1);(1,2)$. We represent, in Figure \ref{fig:CollageMonoidsSN}, the gluing condition
    $$ \Gamma_{\sigma_1} + M(\tau_{1,2},\Gamma_{\sigma_1}) = \{0\}\times 2\Z + \N_{>0}\times \Z = \Gamma_{\sigma_2} + M(\tau_{1,2},\Gamma_{\sigma_2}) $$
    and one can see on the picture that the condition is verified because all the faces are saturated and because $M(\tau_{1,2},\Gamma_{\sigma_2})$ is a subgroup of finite index of the group associated to $\Int(\Gamma_{\sigma_2})$.
    
    \begin{figure}[h]
    \centering
    \begin{tabular}{cc|cc}
    \begin{tikzpicture}[scale=0.7]
    \tikzmath{ integer \nb; \nb = 3;}
    
    \tikzmath{ integer \nbd; \nbd = \nb-1;}
        \foreach \i in {0,...,4}{%
            \foreach \j in {-3,...,4}{%
            \node at (\i, \j) {$\bullet$};
            }
        }

    \node at (1, 4) {{\color{white}$\bullet$}};
    \node at (1, 3) {{\color{white}$\bullet$}};

    \node at (0, -1) {{\color{white}$\bullet$}};
    \node at (0, -1) {$\times$};
    \node at (0, -3) {{\color{white}$\bullet$}};
    \node at (0, -3) {$\times$};
    \node at (0, 1) {{\color{white}$\bullet$}};
    \node at (0, 1) {$\times$};
    \node at (0, 3) {{\color{white}$\bullet$}};
    \node at (0, 3) {$\times$};
    
    \draw[->] (0,0) -- (1,2);
    \draw[dashed, -] (0,0) -- (2,4);
    \draw[->] (0,0) -- (0,-2);
    \draw[dashed, -] (0,0) -- (0,-3.2);

    \end{tikzpicture}

    &&&

    \begin{tikzpicture}[scale=0.7]
    \tikzmath{ integer \nb; \nb = 3;}
    
    \tikzmath{ integer \nbd; \nbd = \nb-1;}
        \foreach \i in {0,...,4}{%
            \foreach \j in {\i,...,4}{%
            \node at (\i, \j) {$\bullet$};
            }
        }

    \node at (0, 1) {{\color{white}$\bullet$}};
    \node at (0, 1) {$\times$};
    \node at (0, 3) {{\color{white}$\bullet$}};
    \node at (0, 3) {$\times$};
    \node at (2, 1) {$\times$};
    \node at (0, -1) {$\times$};
    \node at (0, -3) {$\times$};

    \node at (3, 2) {$\bullet$};
    \node at (4, 2) {$\bullet$};
    \node at (4, 3) {$\bullet$};
    \node at (0, -2) {$\bullet$};
    
    \draw[->] (0,0) -- (0,1.9);
    \draw[dashed, -] (0,0) -- (0,4);
    \draw[->] (0,0) -- (3.95,1.95);

    \end{tikzpicture}\\
    
    $\Gamma_{\sigma_1} + M(\tau_{1,2},\Gamma_{\sigma_1})$ &&& $\Gamma_{\sigma_2} + M(\tau_{1,2},\Gamma_{\sigma_2})$
    \end{tabular}
    
    \caption{Illustration of \Cref{CollageMonoidesSN} }
    \label{fig:CollageMonoidsSN}
\end{figure}

\end{Ex}

%\listoftodos
\bibliographystyle{alpha} 
\bibliography{biblio}
\vspace{1cm}

\end{document}